\documentclass[12pt]{article}

\usepackage{amsmath,amsthm}

\usepackage{amssymb}

\newtheorem{thm}{Theorem}
\newtheorem{prop}[thm]{Proposition}
\newtheorem{lem}[thm]{Lemma}
\newtheorem{cor}[thm]{Corollary}
\newcommand{\FF}{\mathbb{F}}
\newcommand{\ww}{\omega}
\newcommand{\vv}{\bar{\omega}}
\newcommand{\nexteq}{\displaybreak[0]\\ &=}
\newcommand{\cC}{\mathcal{C}}

\DeclareMathOperator{\Aut}{Aut}
\DeclareMathOperator{\wt}{wt}
\DeclareMathOperator{\Nei}{Nei}

\begin{document}
%
\title{Classification of Quaternary Hermitian Self-Dual 
Codes of Length 20}

\author{Masaaki Harada%
\thanks{This work of the first author was supported by JST PRESTO program.}
\thanks{M. Harada is with the 
Department of Mathematical Sciences,
Yamagata University,
Yamagata 990--8560, Japan, and
PRESTO, Japan Science and Technology Agency (JST), Kawaguchi,
Saitama 332--0012, Japan
email: mharada@sci.kj.yamagata-u.ac.jp.}
and Akihiro Munemasa
\thanks{A. Munemasa is with the Graduate School of Information Sciences,
Tohoku University,
Sendai 980--8579, Japan
email: munemasa@math.is.tohoku.ac.jp}
}

\maketitle

\begin{abstract}
A classification of quaternary Hermitian self-dual codes of length
$20$ is given.
Using this classification, a classification of
extremal quaternary Hermitian self-dual codes of length
$22$ is also given.
\end{abstract}


\section{Introduction}
Let
$\FF_4=\{ 0,1,\ww , \vv  \}$ be the finite field of order
four, where $\vv=  \omega^2 = \omega +1$.
Codes over $\FF_4$ are often called quaternary.
All codes in this note are quaternary.
The \textit{Hermitian dual code} $C^{\perp}$ of a code 
$C$ of length $n$ is defined as
$
C^{\perp}=
\{x \in \FF_4^n \mid x \cdot c = 0 \text{ for all } c \in C\},
$
where $x \cdot y = \sum_{i=1}^{n} x_i {y_i}^2$
for $x=(x_1,\ldots,x_n), y=(y_1,\ldots,y_n) \in \FF_4^n$,
which is known as the Hermitian inner product.
A code $C$ is called \textit{Hermitian self-orthogonal} 
if $C \subset C^{\perp}$,
and $C$ is  called \textit{Hermitian self-dual} if $C = C^{\perp}$. 
It is  known that an $[n,k]$ code $C$ is Hermitian self-dual 
if and only if the weights of all codewords of $C$ are even, that is,
$C$ is even, and $n=2k$~\cite{MOSW}. 
It was shown in~\cite{MOSW} that
the minimum weight $d$ of a Hermitian self-dual code of
length $n$ is bounded by
$d \leq 2 \lfloor n/6 \rfloor +2$.
A Hermitian self-dual code of length $n$ and minimum weight 
$d=2 \lfloor n/6 \rfloor +2$ is called \textit{extremal}. 

Two codes $C$ and $C'$ are \textit{equivalent} if there is some monomial
matrix $M$ over $\FF_4$ such that $C' =C M =\{c M \mid c \in C \}$~\cite{MOSW}. 
A monomial matrix which maps $C$ to itself is called an automorphism 
of $C$, and  the set of all automorphisms of $C$ forms the 
automorphism group $\Aut(C)$ of $C$.
Clearly, $\ww I$ and $\vv I$ are elements of $\Aut(C)$,
where $I$ denotes the identity matrix, 
so $\{I,\ww I, \vv I\}$ is the smallest possible
automorphism group.
Such an automorphism group is called {\em trivial}.

The classification of Hermitian self-dual codes was
begun by~\cite{MOSW} and the classification for lengths
up to $14$ was done in~\cite{MOSW}.
The classification is extended to length $18$~\cite{CPS,HLMT}.
At length $20$, a classification of extremal
Hermitian self-dual codes was completed~\cite{Huffman}
under a weaker equivalence (see Subsection~\ref{Subsec:weak}
for the definition).
At length $22$, at least $46$ inequivalent extremal
self-dual codes are known~\cite{Kim}.
In Table~\ref{Tab:C},
the numbers $\#_d$ of inequivalent Hermitian self-dual codes
with minimum weight $d$ are given along with references.

\begin{table}[thb]
\caption{Hermitian self-dual codes}
\label{Tab:C}
\begin{center}
\begin{tabular}{|c|cccc|cc|c|c|}
\hline
&  \multicolumn{4}{c|}{Indecomposable}
& \multicolumn{2}{c|}{Decomposable} 
& \multicolumn{1}{c|}{}& \multicolumn{1}{c|}{}\\
$n$ & $\#_2$& $\#_4$& $\#_6$& $\#_8$&
$\#_2$& $\#_4$& Total &\multicolumn{1}{c|}{References}\\
\hline 
 2 & 1 &    &   &   &    &  &   1 &  \cite{MOSW} \\
 4 & 0 &    &   &   &  1 &  &   1 &  \cite{MOSW} \\
 6 & 0 &  1 &   &   &  1 &  &   2 &  \cite{MOSW} \\
 8 & 0 &  1 &   &   &  2 &  &   3 &  \cite{MOSW} \\
10 & 0 &  2 &   &   &  3 &  &   5 &  \cite{MOSW} \\
12 & 0 &  4 &   &   &  5 & 1&  10 &  \cite{MOSW} \\
14 & 0 &  9 & 1 &   & 10 & 1&  21 &  \cite{MOSW} \\
16 & 0 & 27 & 4 &   & 21 & 3&  55 &  \cite{CPS}  \\
18 & 0 &152 &30 & 1 & 55 & 7& 245 &  \cite{HLMT,Huffman} \\
\hline
20 & 0 &2163&999& 2 &245 &18& 3427& \cite{Huffman}, Section~\ref{Sec:C} \\
22 & 0 &  ? & ? &723&3427&52& ? & Section~\ref{Sec:22} \\
\hline
\end{tabular}
\end{center}
\end{table}

The main aim of this note is to give a complete classification
of Hermitian self-dual codes of length $20$.

\begin{thm}\label{thm:20}
There are $3427$ inequivalent Hermitian self-dual codes of length $20$.
Of these two are extremal, $999$ have minimum weight $6$,
$2181$ have minimum weight $4$, and $245$ have minimum weight $2$.
\end{thm}

The method used,
which is similar to that given in~\cite{Huffman},
is described in Subsection~\ref{Subsec:CM}.

From the classification of Hermitian self-dual $[20,10,6]$
codes,
we also classify extremal Hermitian self-dual 
codes of length $22$.

It has been a question to determine if there is a Hermitian
self-dual code with a trivial automorphism group
(see~\cite[Open Problem (4)]{MOSW}).
There is no Hermitian
self-dual code with a trivial automorphism group
for lengths up to $18$ (see~\cite{MOSW,CPS,HLMT}). 
From our classification given in Section  
\ref{Sec:C},
we see that such a code exists at length $20$.

Generator matrices of all Hermitian self-dual codes of length $20$
and all extremal Hermitian self-dual codes of length $22$
can be obtained electronically from~\cite{Data}.
All computer calculations in this note
were done by {\sc Magma}~\cite{Magma}.

\section{Preliminaries}
\subsection{Classification method}
\label{Subsec:CM}
Here we describe a method for classifying Hermitian self-dual 
codes.
This method is similar to that given in~\cite{Huffman}.

Suppose that 
$C$ is a Hermitian self-dual 
$[n,n/2,d]$ code with $d \ge 4$.
Define a subcode of $C$ as follows
\[
C_0=\{(x_1,x_2,\ldots,x_{n}) \in C \mid x_{n-1}=x_{n}\}.
\]
Since $C^\perp$ has no codeword of weight $2$,
$C_0$ has dimension $n/2-1$.
We may assume that there is a codeword $x=(x_1,\ldots,x_{n})$
of weight $d$ in $C$ with $x_{n-1}=x_{n}\neq0$.
Then, the following code 
\[
C_1=\{(x_1,x_2,\ldots,x_{n-2}) \mid (x_1,x_2,\ldots,x_{n}) \in C_0\}
\]
is a self-dual $[n-2,n/2-1,d-2]$ code.
Thus, the subcode $C_0$ has generator matrix of the form
\begin{equation}\label{eq:G}
G_0=
\left(\begin{array}{ccccc|cc}
&      &      & & &a_1    &a_1   \\
&      &G_1   & & &\vdots &\vdots\\
&      &      & & &a_{n/2-1} &a_{n/2-1}\\
\end{array}\right),
\end{equation}
where $G_1$ is a generator matrix of $C_1$
and $a_i \in \FF_4$ $(i=1,\ldots,n/2-1)$.
It follows that
any Hermitian self-dual $[n,n/2,d]$
code is constructed as
the code $\langle C_0, x \rangle$ 
for some code $C_0$ with generator
matrix of the form (\ref{eq:G})
and some vector $x \in C_0^\perp \setminus C_0$,
where $\langle C_0, x \rangle$ denotes the code
generated by the
codewords of $C_0$ and $x$.

In this way, 
all Hermitian self-dual $[n,n/2,d]$ codes,
which must be checked further for equivalence, 
are constructed, by taking generator matrices of all inequivalent 
self-dual $[n-2,n/2-1,d-2]$ codes
as matrices $G_1$, and by considering $a_1 \in \{0,1\}$ and 
$a_i \in \FF_4$ $(i=2,\ldots,n/2-1)$ in (\ref{eq:G}).
As described in~\cite{Huffman},
the number of possibilities for $a_i$ $(i=2,\ldots,n/2-1)$
is decreased 
by applying elements of
$\Aut(C_1)$ to the first $n-2$ coordinates of (\ref{eq:G}).

\subsection{Mass formula for weight enumerators}
Now we give a mass formula for weight enumerators
of Hermitian self-dual codes.

\begin{lem}\label{lem:wt}
Let $n$ be an even positive integer. 
Let $W_C(y)$ denote the weight enumerator of a code $C$.
Then
\begin{multline}\label{eq:mass} 
\sum_{C}W_C(y)=\\
\prod_{i=0}^{n/2-1}(2^{2i+1}+1)+
\sum_{j=1}^{n/2}
\binom{n}{2j}3^{2j}
\prod_{i=0}^{n/2-2}(2^{2i+1}+1)
y^{2j},
\end{multline} 
where $C$ runs through
the set of all Hermitian self-dual
codes of length $n$.
\end{lem}
\begin{proof}
Let $\wt(x)$ denote the weight of a vector $x \in \FF_4^n$.
\begin{align*}
&\sum_{C}W_C(y)
\\&=
\sum_{j=0}^{n/2} \sum_C\#\{x\in C\mid \wt(x)=2j\} y^{2j}
\nexteq
\sum_{j=0}^{n/2}
\sum_{\substack{x \in \FF_4^n\\ \wt(x)=2j}}\#\{C\mid x\in C\}y^{2j}
\nexteq
\prod_{i=0}^{n/2-1}(2^{2i+1}+1)+
\sum_{j=1}^{n/2}\sum_{\wt(x)=2j}\#\{C\mid x\in C\}y^{2j}
\\&
\text{ (by~\cite[Theorem 19]{MOSW})}
\nexteq
\prod_{i=0}^{n/2-1}(2^{2i+1}+1)+
\sum_{j=1}^{n/2}\sum_{\wt(x)=2j}
\prod_{i=0}^{n/2-2}(2^{2i+1}+1)
y^{2j}     
\\&\text{ (by~\cite[Theorem 23]{MOSW})}
\nexteq
\prod_{i=0}^{n/2-1}(2^{2i+1}+1)+
\sum_{j=1}^{n/2}
\binom{n}{2j}3^{2j}
\prod_{i=0}^{n/2-2}(2^{2i+1}+1)
y^{2j}.
\end{align*}
\end{proof}

\begin{lem}\label{lem:wtC}
Let $n$ and $d$ be even positive integers. Let
$\cC$ be a family of 
pairwise inequivalent
Hermitian self-dual codes of length
$n$ with minimum weight at most $d$. Then $\cC$ is
a complete set of representatives for equivalence classes
of Hermitian self-dual codes of length
$n$ with minimum weight at most $d$, if and only if
\begin{multline}\label{eq:wtd} 
\sum_{C\in\cC}\frac{3^n n!}{\#\Aut(C)}
\#\{x\in C\mid \wt(x)=d\}\\
=\binom{n}{d}3^{d}
\prod_{i=0}^{n/2-2}(2^{2i+1}+1).
\end{multline} 
\end{lem}
\begin{proof}
Consider the coefficient of $y^d$ in the formula 
(\ref{eq:mass}) in Lemma~\ref{lem:wt}.
\end{proof}

\section{Classification of self-dual codes of length 20}
\label{Sec:C}
In this section, we give a classification of 
Hermitian self-dual codes of length $20$.

\subsection{Decomposable codes}
We first consider decomposable Hermitian self-dual codes.
By Theorem 28 in~\cite{MOSW}, any Hermitian self-dual code
with minimum weight $2$ is decomposable as $C_2\oplus C_{18}$,
where $C_2$ is the unique Hermitian self-dual code of length $2$
and $C_{18}$ is some Hermitian self-dual code of length $18$.
Hence, from Table~\ref{Tab:C}, 
there are $245$ inequivalent Hermitian self-dual codes
with minimum weight $2$.
Every decomposable Hermitian self-dual code
with minimum weight $4$ is a direct sum of indecomposable
codes of length at least $6$, since there are no 
Hermitian self-dual code with minimum
weight at least $4$ for lengths less than $6$.
From Table~\ref{Tab:C}, the numbers of indecomposable
Hermitian self-dual codes of lengths $6,8,10,12,14$ are
$1,1,2,4,10$, respectively. It follows that the numbers of
Hermitian self-dual codes of the forms
$C_6 \oplus C_6 \oplus C_8$,
$C_{6} \oplus C_{14}$, 
$C_{8} \oplus C_{12}$ or
$C_{10} \oplus C_{10}$, where
$C_n$ denotes an indecomposable code of length $n$,
are $1,10,4,3$, respectively. Therefore,
there are $18$ inequivalent decomposable Hermitian self-dual codes 
with minimum weight $4$.
There is no decomposable Hermitian self-dual code with minimum
weight $d \ge 6$.

\subsection{Indecomposable codes}

{From} the set of inequivalent Hermitian self-dual $[18,9,2]$ codes
given in~\cite{HLMT}, the method given in Subsection~\ref{Subsec:CM}
allows to enumerate 
the set $\mathcal{C}_{20,4}$ of the $2181$
inequivalent Hermitian self-dual $[20,10,4]$ codes.
Let $\mathcal{C}_{20,2}$ denote the set of the $245$ inequivalent
Hermitian self-dual codes of length $20$ and 
minimum weight $2$. Setting $\mathcal{C}=\mathcal{C}_{20,2}
\cup\mathcal{C}_{20,4}$ in 
Lemma \ref{lem:wtC}, one can verify that
there is no other Hermitian self-dual $[20,10,4]$ code,
by calculating the 
summand in the right-hand side of (\ref{eq:wtd}) for
all codes of $\mathcal{C}$.
Similarly, we found the set of $999$ (resp.\ $2$)
inequivalent Hermitian self-dual codes with minimum weight $6$
(resp.\ $8$).
In this way,
we found the set $\mathcal{C}_{20}$ of $3182$
inequivalent Hermitian self-dual codes 
with minimum weight $d \ge 4$ 
satisfying
\[
\sum_{C \in \mathcal{C}_{20}
\cup \mathcal{C}_{20,2}} 
\frac{3^{20} \cdot 20!}{\#\Aut(C)} 
= 2229034892015508532492061011707,
\]
which is the  constant term of (\ref{eq:mass}).
This constant term $\prod_{i=0}^{9}(2^{2i+1}+1)$
gives the number of the distinct self-dual
codes of length $20$.
The mass formula shows that there is no other Hermitian self-dual code
of length $20$ and the classification is complete.
Therefore, we have Theorem \ref{thm:20}.

This computation was performed in {\sc Magma}~\cite{Magma}.
In principle, such a computation can be done by classifying
Hermitian self-dual codes
by the {\sc Magma} function
{\tt IsIsomorphic}, then their automorphism groups can be 
calculated by {\tt AutomorphismGroup}. 
The orders of the automorphism groups of the $3427$ codes
are listed in Table~\ref{Tab:Aut},
where $(\#\Aut,N(\#\Aut))$ lists the number $N(\#\Aut)$ 
of the codes with an automorphism group of order $\#\Aut$.

Alternatively, the set of $3182$ inequivalent Hermitian self-dual
codes of length $20$ and minimum weight $d\ge4$ can be 
found by a method  similar to the one given in~\cite{HLMT} as 
follows.
Recall that 
two self-dual codes $C$ and $C'$ of length $n$
are said to be {\em neighbors} if the dimension of
$C \cap C'$ is $n/2-1$.
Let $D_{20}$ be the extremal Hermitian self-dual code of length $20$
generated by the second 
generator matrix in~\cite[Fig.~4]{Huffman}.
Let $N_0=\{D_{20}\}$ and
$N_{i+1}= \bigcup_{C \in N_i}\Nei(C)$ $(i=0,1,2,3)$,
where
$\Nei(C)$ denotes the set of inequivalent Hermitian self-dual
neighbors with minimum weight $d \ge 4$ of $C$.
Then the set
\[
N_0 \cup N_1 \cup N_2 \cup N_3 \cup N_4,
\]
contains $3182$ inequivalent Hermitian self-dual
codes of length $20$ and minimum weight $d\ge4$.

\begin{table}
\caption{Orders of the automorphism groups for length $20$}
\label{Tab:Aut}
\begin{center}
{\tiny
\begin{tabular}{|c|l|}
\hline
$d$ & \multicolumn{1}{c|}{$(\#\Aut,N(\#\Aut))$} \\
\hline
8 & (4320, 1), (5760, 1)\\
\hline
6 & 
(3, 419), (6, 328), (9, 13), (12, 103), (15, 1), (18, 26), (24, 15), 
(27, 1), (30, 5), (36, 33), (48, 9), (54, 6), (60, 3), (72, 11), (96, 2), 
\\ &
(108, 4), (144, 3), (192, 3), 
(216, 2), (288, 4), (360, 2), (384, 1), (576, 1), 
(1440, 1), (2160, 1), (3456, 1), 
(5760, 1) 
\\
\hline
4 & 
(12, 377), (24, 258), (36, 22), (48, 363), (72, 38), (96, 134), 
(144, 41), (192, 176), (216, 4), (240, 3), (288, 60), (384, 82), 
(432, 13), 
\\&
(576, 31), (720, 1), (768, 44), (864, 9), (1008, 1), (1080, 2), 
(1152, 92), (1296, 3), (1440, 3), (1536, 30), (1728, 9), 
(2160, 3), 
\\ &
(2304, 42), (2592, 1), (2880, 3), (3072, 6), (3456, 9), (3600, 1), 
(4320, 2), (4608, 31), (5184, 3), (5760, 4), (6144, 1), 
(6480, 1), 
\\&
(6912, 42), (7560, 2), (8640, 1), (9216, 22), (10368, 3), (10800, 1), 
(11520, 3), (12288, 1), (13824, 18), (15360, 1), (16128, 1), 
\\&
(17280, 9), (18432, 3), (20736, 1), (23040, 4), (24576, 1), (25920, 1), 
(27648, 9), (34560, 1), (36288, 1), (36864, 2), (39312, 1), 
\\&
(41472, 8), (48384, 2), (51840, 1), (55296, 21), (69120, 8), (82944, 7), 
(92160, 1), (96768, 1), (103680, 3), (110592, 3), (124416, 1), 
\\&
(138240, 3), (145152, 1), (165888, 4), (184320, 1), (193536, 1), 
(207360, 4), (221184, 4), (248832, 2), (290304, 1), (331776, 3), 
\\&
(345600, 1), (368640, 2), (387072, 1), (414720, 2), (518400, 1), 
(552960, 3), (580608, 2), (622080, 2), (663552, 5), (725760, 1), 
\\&
(777600, 1), (829440, 5), (884736, 2), (1105920, 2), (1161216, 1), 
(1244160, 1), (1382400, 1), (1658880, 1), (1866240, 1), 
\\&
(2073600, 2), (2488320, 1), (3110400, 1), (3538080, 1), (3732480, 1), 
(4064256, 1), (4147200, 2), (4354560, 1), (5806080, 1), 
\\&
(6635520, 1), (7464960, 1), (10886400, 1), (11612160, 1), (13271040, 1), 
(13934592, 1), (14929920, 1), (15552000, 1), 
\\&
(22118400, 1), (23224320, 1), (37324800, 1), (41472000, 1), 
(41803776, 1), (58060800, 1), (66355200, 1), (121927680, 1), 
\\ &
(124416000, 1), (182891520, 1), (278691840, 1), (279936000, 1), 
(933120000, 1), (1045094400, 1), (5573836800, 1), 
\\&
(9405849600, 1)
\\
\hline
2 & 
(36, 1), (72, 5), (108, 1), (144, 9), (162, 1), (216, 6), (288, 4), 
(324, 2), (432, 7), (576, 8), (864, 5), (1080, 2), (1152, 9), (1296, 2), 
\\ &
(1728, 10), (2304, 6), (2592, 1), (3024, 1), (3456, 3), (3888, 1), 
(4608, 2), (5184, 2), (6480, 1), (6912, 13), (7776, 1), (9216, 1), 
\\ &
(10368, 3), (13824, 6), (17496, 1), (18432, 1), (20736, 3), (25920, 1), 
(27648, 4), (31104, 1), (34560, 1), (36864, 1), (41472, 7), 
\\ &
(55296, 1), (62208, 2), 
(72576, 1), (77760, 1), (82944, 8), (84240, 1), 
(103680, 2), (110592, 3), (124416, 2), (139968, 1), 
\\ &
(146880, 1), (165888, 4), (207360, 2), 
(248832, 1), (290304, 1), (311040, 1), (331776, 4), (373248, 1), 
(414720, 1), (497664, 3), 
\\ &
(622080, 1), (746496, 1), (829440, 2), 
(870912, 1), (995328, 3), (1244160, 3), 
(1327104, 3), (1492992, 4), (2239488, 1), 
\\ &
(2488320, 1), (3110400, 1), (3317760, 1), 
(3732480, 1), (3981312, 1), (4245696, 1), 
(4354560, 1), 
(4478976, 1), (5806080, 1),
\\ &
(6635520, 1), (6967296, 1), (8957952, 1), 
(9953280, 2), (11943936, 1), (12192768, 1), 
(15925248, 1), (17915904, 1), (18662400, 1), 
\\ &
(22394880, 1), (24883200, 3), (26127360, 1), 
(37324800, 1), (44789760, 1), 
(67184640, 1), 
(69984000, 1), (107495424, 1), 
\\ &
(139345920, 1), (195955200, 1), 
(219469824, 1), (322486272, 1), (335923200, 1), 
(447897600, 2), (522547200, 1), 
(940584960, 1), 
\\ &
(1114767360, 1), (1254113280, 1), 
(1672151040, 1), (1679616000, 1), 
(2149908480, 1), (2341011456, 1), (5374771200, 1), 
\\ &
(5643509760, 1), 
(20155392000, 1), 
(45349632000, 1), (72559411200, 1), (135444234240, 1), 
(1523747635200, 1), 
\\ &
(219419659468800, 1) 
\\
\hline
\end{tabular}
}
\end{center}
\end{table}
\subsection{A weaker equivalence}
\label{Subsec:weak}
In the above classification, 
we employ monomial matrices over $\FF_4$ in 
the definition for equivalence of codes.
To define a weaker equivalence,
one could consider a conjugation $\gamma$ of $\FF_4$ sending 
each element to its square, in the definition of equivalence,
that is,
two codes $C$ and $C'$ are weakly equivalent if there is some monomial
matrix $M$ over $\FF_4$ such that $C' =C M$
or $C'=CM\gamma$ (see~\cite{Huffman}).
In fact, the classification of extremal self-dual codes of length $20$
in~\cite{Huffman} was done under the weaker equivalence.
Our classification shows that the equivalence classes 
of such codes are the same under both definitions.

We have verified that 
there are $15$, $636$ and $323$ pairs of self-dual codes with minimum
weights $2,4$ and $6$, respectively, under the weaker equivalence.
Hence, there are $3427-(15+636+323)=2453$
self-dual codes under the weaker equivalence
for length $20$.

\section{Classification of extremal self-dual codes of length 22}
\label{Sec:22}
{From} the set of inequivalent Hermitian self-dual 
$[20,10,6]$ codes classified in the previous section,
the method given in Subsection~\ref{Subsec:CM}
allows to enumerate extremal Hermitian 
self-dual codes of length $22$.

For every self-dual $[20,10,6]$ code, 
we have verified that the subcode generated by codewords
of weight $6$ has dimension at least $4$.
Thus, we may assume that the first four rows of
a generator matrix $G_1$ in (\ref{eq:G}) have weight $6$.
This yields that $a_1=1$ and $a_i \in \{1,\ww,\vv\}$
$(i=2,3,4)$.
From the $999$ self-dual $[20,10,6]$ codes,
$723$ inequivalent extremal self-dual codes of
length $22$ are obtained.
Therefore, we have the following:

\begin{prop}\label{Prop:22}
There are $723$ inequivalent extremal Hermitian self-dual codes
of length $22$.
\end{prop}

The orders of the automorphism groups of the $723$ codes
are listed in Table~\ref{Tab:Aut22},
where $(\#\Aut,N(\#\Aut))$ lists the number $N(\#\Aut)$ 
of the codes with an automorphism group of order $\#\Aut$.
The code with an automorphism group of order $1330560$
is equivalent to the code $C_{22,P1}$ in~\cite[Table 2]{Miyabayashi}, 
which is one of the
three inequivalent pure double circulant extremal self-dual codes
of length $22$.
We have verified that the automorphism group is isomorphic to
the direct product of
the Mathieu group $M_{22}$ of degree $22$
and the cyclic group of order $3$.

\begin{table}
\caption{Orders of the automorphism groups for length $22$}
\label{Tab:Aut22}
\begin{center}
{\small
\begin{tabular}{|l|}
\hline
\multicolumn{1}{|c|}{$(\#\Aut,N(\#\Aut))$} \\
\hline
(3, 308), (6, 229), (9, 8), (12, 73), (18, 12), 
(24, 39), (30, 2), (36, 19), 
(48, 1),  
\\
(60, 1), (66, 2), (72, 5), (108, 1), 
(120, 1),  (180, 3), (192, 1), 
(240, 1), (288, 4), \\
(324, 1), (360, 1), (384, 2), 
(504, 2), (864, 2), 
(1728, 2), (17280, 2), (1330560, 1) \\
\hline
\end{tabular}
}
\end{center}
\end{table}

We have verified that 
there are $301$ pairs of extremal self-dual codes of
length $22$, under the weaker equivalence given in 
Subsection \ref{Subsec:weak}.
Hence, there are $723-301=422$ extremal self-dual codes of
length $22$, under the weaker equivalence.
We have also verified that the $422$ extremal self-dual codes
have different numbers $(B_0,B_1,\ldots,B_{22})$,
where $B_j$ denotes the number of distinct cosets of weight $j$.
This shows that the $422$  extremal self-dual codes
are certainly inequivalent.

The smallest possible automorphism group is of order $3$.
{From} the constant term of (\ref{eq:mass}),
the number of inequivalent self-dual codes of length $22$
is at least 
\[
\frac{\prod_{i=0}^{10} (2^{2i+1} + 1)}{22! \cdot 3^{21}} 
> 397588.
\]

\section{Trivial automorphism groups}
\label{Sec:Aut}

It has been a question to determine if there is a Hermitian
self-dual code with a trivial automorphism group
(see~\cite[Open Problem (4)]{MOSW}).
There is no Hermitian
self-dual code with a trivial automorphism group
for lengths up to $18$ (see~\cite{MOSW,CPS,HLMT}). 
From our classification given in Section 
\ref{Sec:C},
we see that such a code exists at length $20$ (see 
Table~\ref{Tab:Aut}).
Hence, we have an answer to the above problem.

\begin{cor}\label{cor}
The smallest length for which there is a
quaternary Hermitian self-dual code with a trivial
automorphism group is $20$.
\end{cor}

We describe here an example of a code satisfying the conditions
of Corollary~\ref{cor}.
Let $B_{20}$ be the pure double circulant code with
generator matrix $(\ I\ , \ R \ )$, where
$R$ is the circulant matrix with first row
$(\vv,0,\ww,0,\ww,0,\vv,1,\ww,1)$.
The following code
\[
C_{20}= \langle B_{20} \cap 
\langle v \rangle^\perp, v \rangle
\]
where
$
v=(\ww,1,1,1,1,1,1,1,1,1,0,\vv,0,0,1,0,\ww,\vv,1,\vv),
$
is a Hermitian self-dual code of length $20$ with a trivial
automorphism group.



In addition, there is no extremal Hermitian self-dual code 
with a trivial automorphism group for length $20$, but
there are $308$ extremal Hermitian self-dual codes
with trivial automorphism groups for length $22$
(see Tables~\ref{Tab:Aut}
and \ref{Tab:Aut22}).
Hence, the smallest length for which there is a
quaternary extremal Hermitian self-dual code with a trivial
automorphism group is $22$.
For example, the code with generator matrix $(\ I\ ,\ M\ )$
is an extremal self-dual $[22,11,8]$ code with a trivial 
automorphism group,
where $M$ is given in Figure~\ref{Fig}.



\begin{figure}[!t]
\centering
\[
M=
\left(\begin{array}{ccccccccccc}
  1 &   0 &   1 &   0 &   1 &   0 & \ww & \vv & \ww & \vv &   0 \\
\ww & \ww &   0 & \ww &   1 &   0 & \ww & \ww & \vv & \ww &   1 \\
\vv & \ww & \vv &   0 &   0 & \vv & \ww & \vv & \vv &   1 & \ww \\
\vv &   0 &   1 &   1 &   0 & \ww & \vv &   1 & \ww & \ww & \ww \\
\ww & \ww &   0 & \ww & \vv &   1 & \vv & \vv & \vv &   0 & \vv \\
  1 & \vv &   1 &   0 &   0 & \ww &   0 & \vv &   0 & \vv &   1 \\
  0 &   1 & \vv &   0 &   0 & \ww &   1 &   0 &   1 & \ww & \vv \\
\ww &   0 & \vv & \vv & \vv &   1 & \vv &   0 & \ww & \vv & \vv \\
  1 & \ww & \ww & \ww &   1 &   1 &   0 & \vv &   1 &   0 & \ww \\
  1 & \ww & \vv & \ww & \vv &   1 &   1 & \vv &   0 & \ww &   0 \\
  0 &   0 & \ww & \vv & \vv & \ww &   1 &   1 &   1 & \ww & \ww 
\end{array}\right)
\]
\caption{An extremal self-dual $[22,11,8]$ code with a trivial automorphism group}
\label{Fig}
\end{figure}

\end{document}